\documentstyle[12pt,amsfonts,amscd]{amsart}
\newtheorem{Theorem}{Theorem}[section]
\newtheorem{Definition}[Theorem]{Definition}
\newtheorem{Proposition}[Theorem]{Proposition}

\newtheorem{Lemma}[Theorem]{Lemma}
\newtheorem{Corollary}[Theorem]{Corollary}
\theoremstyle{remark}
\newtheorem{Example}[Theorem]{Example}

\def\CC{{\Bbb C}}
\def\RR{{\Bbb R}}

\def\be{\begin{enumerate}}
\def\ee{\end{enumerate}}
\def\bT{\begin{Theorem}}
\def\eT{\end{Theorem}}
\def\bP{\begin{Proposition}}
\def\eP{\end{Proposition}}
\def\bD{\begin{Definition}}
\def\eD{\end{Definition}}
\def\bE{\begin{Example}}
\def\eE{\end{Example}}
\def\bL{\begin{Lemma}}
\def\eL{\end{Lemma}}
\def\bC{\begin{Corollary}}
\def\eC{\end{Corollary}}

\begin{document}
\title{Testing holomorphy on curves}
\author{Buma L. Fridman and  Daowei Ma}
\begin{abstract} For a domain $D\subset {\Bbb{C}}^n$ we construct a continuous foliation of $D$ into 
one real dimensional curves such that any function $f\in {C^1(D)}$ which can be extended holomorphically into some neighborhood of each curve in the foliation will be holomorphic on $D$.
\end{abstract}


\address{ buma.fridman@@wichita.edu, Department of Mathematics,
Wichita State University, Wichita, KS 67260-0033, USA}
\address{ dma@@math.wichita.edu, Department of Mathematics,
Wichita State University, Wichita, KS 67260-0033, USA}
\maketitle \setcounter{section}{1}

This paper complements the study of the following general question. Let $f$
be a function on a domain $D$ in complex $n$-dimensional space, and its
restrictions on each element of a given family of subsets of $D$ is
holomorphic. When can one claim that $f$ has to be holomorphic in $D$?

This is a natural question arising from the fundamental Hartogs theorem
stating that a function $f$ in ${\Bbb{C}}^n$, $n>1$, is holomorphic if it is
holomorphic in each variable separately, that is, $f$ is holomorphic in ${\Bbb{C}}^n$ if it is holomorphic on every complex line parallel to an axis. The complex lines parallel to an axis form a continuous foliation of ${\Bbb{C}}^n$ into two real dimensional planes. So if a function is holomorphic along each component of these $n$ foliations, then it is holomorphic on ${\Bbb{C}}^n$. We are interested in finding a one family of one real dimensional curves forming a foliation such that a similar theorem will hold. 
There is a body of interesting work on testing the holomorphy property on curves: see [A1-A3, AG, E, G1-G3, T1, T2] and references in those articles. Some of these results assume a holomorphic extension into the inside of each closed curve in a given family, others a ``Morera-type'' property. 

Below we use the following definition. Let $S\subset {\Bbb{C}}^n$. We say that $f:S\rightarrow \Bbb{C}$ is
holomorphic if $f$ is a restriction on $S$ of a function holomorphic in some
open neighborhood of $S$. We prove the following

\bT \label{T1} Let $D\subset {\Bbb{C}}^n$ be a domain. Then there exists a continuous foliation $E$ of $D$ into one (real) dimensional curves, such that any $C^1$ function on $D$ which is holomorphic on each of the curves of $E$, is holomorphic on $D$.
\eT
The needed foliation will be constructed as a homeomorphic image of the natural continuous foliation of $\overline{D}$ by segments parallel to the axis~$\ {Re}z_1$.

First we need the following notion (see [FM]). Let $S\subset {\Bbb{C}}$, $p\in S$. A point $t$ in $T:=\{z\in\CC: |z|=1\}$ is said to be a limit direction of $S$ at $p$ if there exists a sequence $(q_j)$ in $S$ such that $\lim_j q_j=p$ and $\lim_j\tau(p,q_j)=t$, where $\tau(p,q_j):=(q_j-p)/|q_j-p|$.

\bL \label{TS}
Let $U \subset \Bbb{C}$ be an open set, $p\in U \cap S$ and there are at least two limit directions $t_1,t_2$ of $S$ at $p$. Suppose a function $f\in C^1 (U)$ is
holomorphic on $S \cap U$. If $t_1\neq \pm t_2$ then $\frac{\partial f}{\partial \overline{z}}=0$ at $p$.
\eL
\begin{pf}
The derivatives of $f$ along linearly independent directions $t_1$ and $t_2$ coincide with derivatives of a holomorphic function in the neighborhood of $p$. The statement now follows from the Cauchy-Riemann equations. 
\end{pf}
\bE Consider a set $\gamma \subset \Bbb{C}$, which is an angle ($\gamma =\angle$) in a neighborhood of a point $p \in \gamma$  formed by two linear segments. If this angle $\theta$ satisfies $0<\theta <\pi $, then at the tip of the angle $p\in \gamma$, $\gamma$ has two linearly independent directions. 
\eE
In general if in a neighborhood of a point $p \in \gamma$ the curve $\gamma \subset \RR^m$ lies in a two real dimensional plane $M$ and forms an angle $0<\theta <\pi $ there, we will say that $\gamma$ has an angular point at $p$.

For the construction of the continuous foliation in the Theorem~\ref{T1} we also need the following general statement.
\bL \label{per} Let $M$ be a two-dimensional plane in $\RR^m$ with $m\ge2$, let $p$ be a point in $M$, let $U$ be a neighborhood of $p$ in $\RR^m$,  and let $\gamma$ be a $C^\infty$ curve passing through $p$ and relatively closed in $U$. Then there is a homeomorphism $\Phi: \RR^m\to\RR^m$ such that (a) the function $\Phi$ is a ($C^\infty$) diffeomorphism on $\RR^m-\{p\}$, (b) the restriction of $\Phi$ to $\RR^m-U$ is the identity map, (c) a neighborhood of $\Phi(p)$ in $\Phi(\gamma)$ lies in $M$, and (d) the curve $\Phi(\gamma)$ has an angular point at $\Phi(p)$.
\eL

\begin{pf} Let $e_1, e_2, \dots, e_m$ be the standard basis of $\RR^m$, {\it i.e.}, $e_1=(1,0,\dots,0)$, $e_2=(0,1,0,\dots, 0)$, etc. Choose a vector $v$ parallel to $M$ such that $v$ and the tangent vector of $\gamma$ at $p$ are linearly independent. Without loss of generality, we assume that $p=0$, $v=e_2$, and $M$ is spanned by $e_1$ and $e_2$. Let $S=\RR e_2=\{te_2\in\RR^m: t\in\RR\}$. For $r>0$ let $B_r$ denote the open ball in $\RR^m$ of center $0$ and radius $r$. There is a neighborhood $V$ of $0$, a $\delta>0$, and a diffeomorphism $G: V\to B_{3\delta}$ such that $V\subset\subset U$, $G(\gamma\cap V)=\RR e_1\cap B_{3\delta}$, and $G(S\cap V)=\RR e_2\cap B_{3\delta}$.

Let $\omega: \RR\to\RR$ be a $C^\infty$ function such that 
$0\le \omega(t)\le 1$ for all $t$, $\omega(t)=0$ for $|t|\ge 2$, and $\omega(t)=1$ for $|t|\le 1$. 
Define a vector field $X$ on $\RR^m$ by $X(y)=\omega(|y|/\delta)(y_1+y_2)(e_2-e_1)$, where $y=\sum_{j=1}^m  y_je_j$. Let $\theta: \RR\times\RR^m\to\RR^m$ be the associated action. Define $\lambda: \RR^m\to\RR^m$ by $\lambda(y)=\theta(1,y)$. Then $\lambda$ is a diffeomorphism, and the restriction of $\lambda$ to $\RR^m-B_{2\delta}$ is the identity map, since $X=0$ there.
We claim that for $-\delta <s<\delta$, $\lambda(se_1)=se_2$. Indeed, it is straightforward to verify that the curve $\tau(t)=s(1-t)e_1+ste_2$ satisfies $\tau(0)=se_1$, $\tau(1)=se_2$, and $\tau'(t)=X(\tau(t))$ for $0\le t\le1$, and hence $\tau|_{[0,1]}$ is a segment of an integral curve of $X$. 

Define a diffeomorphism $g: \RR^m\to\RR^m$ by

\begin{equation*} 
   g(x)= \left\{
        \begin{array}{ll}
        G^{-1}\circ\lambda\circ G(x),&\mbox{if $x\in G^{-1}(B_{2\delta})$,}\\
            x,&\mbox{if $x\not\in G^{-1}(B_{2\delta})$.}
        \end{array}
    \right.
\end{equation*}
Then $g(0)=0$, and $V_1\cap g(\gamma)\subset S\subset M$, where $V_1:=G^{-1}(B_\delta)$. 

Choose a $K>0$ so that the function $\psi(t):=K\omega(2t)(1-|t|)$ satisfies the Lipschitz condition $|\psi(t_1)-\psi(t_2)|\le |t_1-t_2|/2$. It is clear that for each $\eta>0$ the function $\psi_\eta(t):= \eta \psi(t/\eta)$ satisfies the same Lipschitz condition. 

Choose an $\eta>0$ such that $B_{\eta}\subset\subset V_1$. Define $h: \RR^m\to\RR^m$ by
$h(x)=x+\psi_\eta(|x|)e_1$. Then $h$ is a homeomorphism, and it is a diffeomorphism away from the origin. For $x\in \RR^m-B_\eta$, $h(x)=x$. The set $h(g(\gamma)\cap B_{\eta/2})$ lies in $M$ and equals
$\{K(\eta-|t|)e_1+te_2: -\eta/2<t<\eta/2\},$
which is the union of two line segments forming an angle $2\tan^{-1}(1/K)$ at the point $K\eta e_1$.

Let $\Phi=h\circ g$. Then $\Phi$ has all the prescribed properties.
\end{pf}

We now proceed with the construction of $E$ and proof of the Theorem~\ref{T1}.

\begin{pf}
Consider $E_0$ a natural continuous foliation of $D$ by segments parallel to the axis $\ {Re}z_1$.

1. Pick a sequence $\{w_k\}\subset D$, such that $\overline{\{w_m\}}_{m\equiv l\pmod{n}}=\overline{D}$ for every $l=1,...,n$. We also can choose the sequence in such a way that no
two points lie on the same line segment of $E_0$, so we assume that each of these points $w_k$ lies on a
unique segment $L_k$. 

2. We now proceed by induction on $k.$

(1). For $k=1$ pick $\varepsilon _1>0$ so small that the ball $\overline{B}
(w_1,\varepsilon _1)\subset D.$ Use Lemma~\ref{per} to create  a homeomorphism $\Phi _1:\overline{D}
\rightarrow \overline{D}$ which is a diffeomorphism on $D\backslash \{w_1\}$ with the following properties. At the point  $\Phi _1(w_1)$ the image  $\Phi _1(L_1)$
has an angle $0<\alpha _1<\pi$, and that angle (as a portion of $\Phi _1(L_1)$) lies in the plane parallel to $z_1$.
Let $d_1={\min_{{d(z,w)\geq 1/2} }}d(\Phi _1(z),\Phi _1(w))$, where $d$ is the Euclidean distance between two
points in ${\Bbb{C}}^n$.

(2). Consider now step $k=s+1$. By now we have constructed a homeomorphism $\Phi _j:
\overline{D}\rightarrow \overline{D}$ which is diffeomorphic on $D\backslash \{w_1,...,w_j\}$, $\varepsilon _j>0$, and $d_j={\min_{{d(z,w)\geq 1/(j+1)} }} d(\Phi _j(z),\Phi _j(w))$
for all $j\leq s$. Also for all $j\leq s$ we assume that $\Phi _s(L_j)=\Phi _j(L_j)$, and that $\Phi _j(L_j)$ has
an angular point at $\Phi _j(w_j)$ in the plane parallel to $z_l$, axis, where $l\equiv j \pmod{n}$. 

Pick now $\varepsilon _{s+1}>0$ such that the following
four conditions hold:

(a) $\varepsilon _{s+1}<\frac 12\varepsilon _s$.

(b)$B(\Phi _s(w_{s+1}),\varepsilon _{s+1})\subset D\backslash \{\cup _{j\leq s}\Phi _s(L_j)\}$.

(c) $\varepsilon _{s+1}<\frac 1{16}d_s$.

(d) $\varepsilon _{s+1}$ is so small that one can use Lemma~\ref{per} to create the
specific perturbation $\widetilde{\Phi }_{s+1}:\overline{D}\rightarrow 
\overline{D}$ inside $B(\Phi _s(w_{s+1}),\varepsilon _{s+1})$ that makes an angle $
0<\alpha _{s+1}<\pi$ at the point  $\widetilde{\Phi }_{s+1}(\Phi _s(w_{s+1}))$ in
the plane parallel to $z_l$ axis, where $l\equiv s+1 \pmod{n}$, and is the identity map 
outside $B(\Phi _s(w_{s+1}),\varepsilon _{s+1})$.

Consider now  $\Phi _{s+1}:\overline{D}\rightarrow \overline{D}$ , which is
defined the following way: $\Phi _{s+1}=\widetilde{\Phi }_{s+1}\circ \Phi _s$.

One can see that $\Phi _{s+1}(z)$ is a well defined homeomorphism which is
diffeomorphic on $D\backslash \{w_1,w_2,...,w_{s+1}\}$. Also for all $j\leq s+1$, $\Phi _{s+1}(L_j)=\Phi _j(L_j)$.
Let $d_{s+1}={\min_{{d(z,w)\geq 1/(s+2)} }}d(\Phi _{s+1}(z),\Phi _{s+1}(w))$.

Consider now $\Phi _0={\lim_j }\Phi _j$. The
limit exists since $\Vert \Phi _{s+1}-\Phi _s\Vert <\frac 1{2^{s-1}}\varepsilon
_1$ for all $s$. We shall prove that $\Phi _0$ is a homeomorphism from $\overline{D}$ onto $
\overline{D}$. All we need to check is that for two points $z\neq w$ in $D$, 
$\Phi _0(z)\neq \Phi _0(w)$. Indeed, find the smallest $s$, such that $
d(z,w) \ge \frac 1{s+1}$. By the construction $d(\Phi _s(z),\Phi _s(w))\geq d_s$
, $\Vert \Phi _{j+1}-\Phi _j\Vert \leq 2 \varepsilon _{j+1}$ for all $j;$ considering the last inequality for $j\geq s$, we have $
d(\Phi _s(z),\Phi _0(z))\leq {2\sum_{j\geq {s+1}} }\varepsilon _j<
{\sum_{j\geq 0} \frac 1{2^{j-1}}\varepsilon _{s+1} = 4 \varepsilon_{s+1}}< \frac 14d_s$. Same inequality holds for the point $w$. So, $d(\Phi _0(z),\Phi
_0(w))>\frac 12d_s>0$, and therefore $\Phi _0(z)\neq \Phi _0(w)$.

We now check that the continuous foliation $E=\Phi _0(E_0)$ satisfies the
theorem. First we notice that for all $j$ by construction $\Phi _0(L_j)=\Phi _j(L_j)$, and therefore $\Phi _0(L_j)$ has
an angular point at $\Phi _0(w_j)$ in the plane parallel to $z_l$ axis, where $l\equiv j \pmod{n}$.

If a function  $f\in C^1(D)$ is holomorphic on each of the curves in $E$, then by Lemma~\ref{TS}, $\frac{\partial f}{\partial \overline{z}_l}=0$
at an everywhere dense set in $D$, and therefore on all of $D$, and for each $
l$. By Hartogs theorem, $f$ is holomorphic on $D$.  \end{pf}

\end{document}